\documentclass[11pt]{article}
\usepackage{amssymb,amsmath,latexsym, verbatim}
\allowdisplaybreaks

\begin{document}

\begin{doublespace}

\newtheorem{thm}{Theorem}[section]
\newtheorem{lemma}[thm]{Lemma}
\newtheorem{defn}[thm]{Definition}
\newtheorem{prop}[thm]{Proposition}
\newtheorem{corollary}[thm]{Corollary}
\newtheorem{remark}[thm]{Remark}
\newtheorem{example}[thm]{Example}
\newtheorem{cond}[thm]{Condition}
\numberwithin{equation}{section}
\def\ee{\varepsilon}
\def\qed{{\hfill $\Box$ \bigskip}}
\def\NN{{\cal N}}
\def\AA{{\cal A}}
\def\MM{{\cal M}}
\def\BB{{\cal B}}
\def\CC{{\cal C}}
\def\LL{{\cal L}}
\def\DD{{\cal D}}
\def\FF{{\cal F}}
\def\EE{{\cal E}}
\def\QQ{{\cal Q}}
\def\RR{{\mathbb R}}
\def\R{{\mathbb R}}
\def\L{{\bf L}}
\def\K{{\bf K}}
\def\S{{\bf S}}
\def\A{{\bf A}}
\def\E{{\mathbb E}}
\def\F{{\bf F}}
\def\P{{\mathbb P}}
\def\N{{\mathbb N}}
\def\eps{\varepsilon}
\def\wh{\widehat}
\def\wt{\widetilde}
\def\pf{\noindent{\bf Proof.} }
\def\beq{\begin{equation}}
\def\eeq{\end{equation}}
\def\lam{\lambda}
\def\ds{\displaystyle}

\title{\Large \bf Fatou and relative Fatou theorem for subordinate Brownian motions with Gaussian components on smooth domains}
\author{Yunju Lee, Hyunchul Park}

\date{\today }
\maketitle

\begin{abstract}
We prove relative Fatou's theorem for nonnegative harmonic functions with respect to a large class of killed subordinate Brownian motions
with Gaussian components in bounded $C^{1,1}$ open sets in $\mathbb{R}^{d}$, $d\geq 2$, which asserts the existence of nontangential limit of the ratio of two harmonic functions with
respect to the killed processes.
When $D=B(x_{0},r)$ is a ball we prove Fatou theorem. That is, we establish the existence of nontangential limit of
a single nonnegative harmonic function. We also prove this is the best result possible by showing that there is a nonnegative harmonic function which does not
have a tangential limit a.e. when $d=2$ and $D=B(0,1)$.
\end{abstract}

\section{Introduction}
In 1906, Fatou proved in \cite{Fatou} that bounded harmonic functions in the open unit disk have nontangential limits almost everywhere on the unit circle. Later Fatou's theorem is extended to a more general setting, namely, relative Fatou theorem.
In \cite{Doob} the author proved that the ratio $u/v$ of two positive harmonic functions defined on an open solid sphere with respect to Brownian motions has nontangential limits almost everywhere with respect to the Martin-representing measure of $v$.

However, when underlying processes are jump processes, Fatou and relative Fatou theorem are not true anymore in the form stated above (see
\cite{BY} for the counterexample for Fatou theorem when the underlying processes are symmetric stable processes).
When the underlying processes are jump processes, relative Fatou theorem means the existence of nontangential limit of the ratio $u/v$ of two nonnegative functions where $v$ is harmonic with respect to killed processes $X^{D}$ and $u$ is harmonic with respect to either
$X$ or $X^{D}$.
For symmetric stable processes, relative Fatou theorem in this form has been established in \cite{BD, MR} for bounded $C^{1,1}$ domains and Lipschitz domains when both $u$ and $v$ are harmonic with respect to killed symmetric stable processes
and in \cite{Kim} for bounded $\kappa$-fat open sets when $u$ is harmonic with respect to symmetric stable processes
and $v$ is harmonic with respect to killed symmetric stable processes, respectively.
Recently relative Fatou theorem is established for a large class of pure jump subordinate Brownian motions in bounded $\kappa$-fat open sets in
\cite{KL}.

In many situations, such as finance and control theory, one needs Markov processes that have both diffusion and jump components.
In an analytical point of view, one would like to study operators that have both local and non-local parts and these operators correspond to processes with both diffusion and jump components.
The prototype of these processes would be an independent sum of Brownian motions and symmetric stable processes, and they are studied in \cite{CKS3, CKS1, CKSV1, CKS2}. The potential theory of a large class of subordinate Brownian motions with Gaussian components is also studied in \cite{KSV1}.
The purpose of this paper is to study the boundary behavior of nonnegative harmonic functions with respect to subordinate Brownian motions that have both diffusion
and jump parts. More precisely we establish the existence of the nontangential limit of the ratio
of nonnegative harmonic functions with respect to a large class of killed subordinate Brownian motions with Gaussian components in bounded $C^{1,1}$ open sets.
When the $C^{1,1}$ open set is a ball, we prove Fatou theorem, which asserts that nonnegative harmonic functions with respect to
the killed subordinate Brownian motions on a ball have nontangential limits almost everywhere on the boundary of the ball.

In this paper we use the recent results about the identification of Martin boundary with Euclidean boundary and the sharp two-sided Green function
estimates for subordinate Brownian motions with Gaussian components in \cite{KSV1} to establish relative Fatou theorem on bounded $C^{1,1}$ open sets.
We follows arguments in \cite{MR} closely while making some modifications if necessary to prove relative Fatou theorem.
The main contribution of this paper is to realize that when the domain is a ball and the processes have a diffusion part,
one can use relative Fatou theorem to establish Fatou theorem.
We investigate the case when the normalizing function corresponds to the Martin integral with respect to the surface measure of the boundary of the domain $D$.
Since $X$ have a diffusion part, it is natural to expect that the probability
that the processes exit $D$ through its boundary is close to 1 as the starting point of the processes approaches $\partial D$.
The function defined by $\P_{x}\left(X_{\tau_{D}}\in \partial D\right)$ is a nonnegative and nonzero harmonic function
with respect to $X^{D}$ and it converges to 1 as $x\rightarrow \partial D$.
In case of $D$ to be a ball $B=B(x_{0},r)$ centered at $x_{0}$, the harmonic measure restricted to $\partial D$ is a normalized surface measure of $\partial D$.
Hence with the normalizing function to be $\P_{x}\left(X_{\tau_{B}}\in \partial B\right)$, relative Fatou theorem on balls can be used to prove
the existence of the nontangential limits of a nonnegative harmonic function with respect to $X^{B}$
whose Martin measure is absolutely continuous with respect to the surface measure of $\partial B$.
We also prove this is the best result possible in a sense that there exists a nonnegative harmonic function
$u(x)$ on the unit ball which fails to have tangential limits for a.e. on $\partial B$.
The investigation of Fatou theorem and harmonic measure of subordinate Brownian motions with Gaussian component in $C^{1,1}$ open sets will appear
in a forthcoming paper.

The rest of the paper is organized as follows. In Section \ref{preliminary} we recall some basic facts about
a large class of subordinate Brownian motions with Gaussian components studied in \cite{KSV1}
and also provide the sharp two-sided Martin kernel estimate
in bounded $C^{1,1}$ open sets $D$ (see \eqref{equn:Martin}).
In Section \ref{RFT} we prove relative Fatou theorem for nonnegative harmonic functions
with respect to these killed subordinate Brownian motions with Gaussian components on bounded $C^{1,1}$ open sets $D$.
In the end of Section \ref{RFT} we consider the most natural case when the function $v(x)$ corresponds to the Martin integral with
respect to the surface measure of the boundary of the domain $D$.
In Section \ref{ball} we investigate the case when $D$ corresponds to a ball $B(x_{0},r)$ centered at $x_{0}$. In this section we establish Fatou theorem
for a ball and shows that Stolz open sets are best possible for Fatou theorem being hold true
by showing that there is a nonnegative and bounded harmonic function whose radial and tangential limits
do not agree for almost every point on $\partial B(x_{0},r)$.

Throughout this paper, we will use $c$ and $c_{i}$ to denote positive constants depending
(unless otherwise explicitly stated)
only on the dimension $d$ but whose value may change from one appearance to another, even within a single line.
We will use $\delta_{D}(x)$ to denote the Euclidean distance between $x$ and $D^{c}$.
For any two positive functions $f$ and $g$, $f\asymp g$ means that there exists a positive constant $c\geq 1$ such that
$c^{-1}g\leq f\leq cg$ on their common domain of the definition.

\section{Preliminaries}\label{preliminary}
Throughout this paper we will assume $d\geq 2$.
In this section we define subordinate Brownian motions with Gaussian components and state some properties about them.
Recall that an one-dimensional L\'evy process $S=(S_{t}, t\geq 0)$ is called a subordinator if it is nonnegative. A subordinator $S$ can be characterized by its Laplace exponent $\phi$ through the relation
$$
\E[e^{-\lam S_{t}}]=e^{-t\phi(\lam)}, \quad t>0, \lam>0.
$$
A smooth function $\phi : (0, \infty) \to [0, \infty)$ is called a Bernstein function if $(-1)^n D^n \phi \le 0$ for every positive integer $n$.
The Laplace exponent $\phi$ of a subordinator is a Bernstein function with $\phi(0+)=0$ and can be written as
$$
\phi(\lam)=b\lam+\int_{(0,\infty)}(1-e^{-\lam t})\mu(dt), \quad \lam >0,
$$
where $b\geq 0$ and $\mu$ is a measure on $(0,\infty)$ satisfying $\int_{(0,\infty)}(1\wedge t)\mu(dt)<\infty$. $\mu$ is called the L\'evy measure of $\phi$. In this paper we will assume that $b>0$ in order to have a nontrivial diffusion part for subordinate Brownian motions.
Without lose of generality we assume $b=1$.

Suppose that $W=(W_{t}: t\geq 0)$ is a $d$-dimensional Brownian motion and $S=(S_{t}: t\geq 0)$ is a subordinator
with Laplace exponent $\phi$, which is independent of $W$.
The process $X=(X_{t}:t\geq 0)$ defined by $X_{t}=W(S_{t})$ is called a subordinate Brownian motion and its infinitesimal generator is given by $\phi(\Delta):=-\phi(-\Delta)$, which can be constructed via Bochner's functional calculus. On $C_{b}^{2}(\R^{d})$ (the collection of $C^{2}$
functions in $\R^{d}$ which, along with partial derivatives up to order 2, are bounded),
$\phi(\Delta)$ is an integro-differential operator of the type
$$
\Delta f(x)+\int_{\R^{d}}\left(f(x+y)-f(x)-\nabla f(x)\cdot y 1_{\{|y|\leq 1\}}\right)J(dy),
$$
where the measure $J$ has the form $J(dy)=j(|y|)dy$ with $j:(0,\infty)\rightarrow (0,\infty)$ given by
$$
j(r)=\int_{0}^{\infty}(4\pi t)^{-d/2}e^{-r^{2}/(4t)}\mu(dt).
$$

Throughout this paper we will impose two conditions on $\phi$ and $\mu$.
\begin{cond}\label{con:con}
\begin{enumerate}
\item
The Laplace exponent $\phi$ of $S$ is a completely Bernstein function. That is, the L\'evy measure $\mu$ has a completely monotone density
(i.e., $\mu (dt) = \mu(t)dt$ and $(-1)^n D^n \mu \ge 0$ for every non-negative integer $n$).
\item For any $K>0$, there exists $c=c(K)>1$ such that
\beq\label{cond:2}
\mu(r)\leq c\mu(2r)
\quad \text{for }\,\,  r\in (0,K).
\eeq
\end{enumerate}
\end{cond}
Note that Condition \ref{con:con} is the main assumption imposed in \cite{KSV1}.

For any open set $D\subset \R^{d}$, $\tau_{D}:=\inf\{t>0 : X_{t}\notin D\}$ denotes the first exit time from $D$ by $X$.
We will use $X^{D}$ to denote the process defined by $X_{t}^{D}(\omega)=X_{t}(\omega)$
if $t<\tau_{D}(\omega)$ and
$X_{t}^{D}(\omega)=\partial$ if $t\geq \tau_{D}(\omega)$, where $\partial$ is a cemetery point.
It is well known that $X^{D}$ is a strong Markov process with state space $D\cup \{\partial\}$.
For any function $u(x)$ defined on $D$ we extend it to $D\cup \{\partial\}$ by letting $u(\partial)=0$.
It follows from \cite[Chapter 6]{BBKRSV} that the process $X$ has a transition density $p(t,x,y)$ which is jointly continuous. Using this and the strong Markov property, one can easily check that
$$
p_{D}(t,x,y):=p(t,x,y)-\E_{x}[p(t-\tau_{D},X_{\tau_{D}},y); t>\tau_{D}], \quad x,y\in D
$$
is continuous and is a transition density of $X^{D}$.
For any bounded open set $D\subset \R^{d}$, we will use $G_{D}(x,y)$ to denote the Green function of $X^{D}$, i.e.,
$$
G_{D}(x,y):=\int_{0}^{\infty}p_{D}(t,x,y)dt, \quad x,y\in D.
$$
Note that $G_{D}(x,y)$ is continuous on $\{(x,y)\in D\times D : x\neq y\}$.

The L\'evy density is given by $J(x,y)=j(|x-y|)$, $x,y\in \R^{d}$ and it determines a L\'evy system for $X$, which describes the jumps of the process $X$:
For any nonnegative measurable function $f$ on $\R_{+}\times\R^{d}\times \R^{d}$ with $f(s,x,x)=0$ for all $s>0$ and $x\in\R^{d}$, and stopping time $T$
with respect to $\{{\cal F}_t: t\ge 0\}$,
$$
\E_{x}\left[\sum_{s\leq T}f(s,X_{s-},X_{s})\right]=\E_{x}\left[\int_{0}^{T}\left(\int_{\R^{d}}f(s,X_{s},y)J(X_{s},y)dy\right)ds\right].
$$
Using L\'evy system, we know that for any nonnegative function $f\geq 0$ and every bounded open set $D$ we have
\beq\label{eqn:Levy system}
\E_{x}\left[f(X_{\tau_{D}}), X_{\tau_{D^{-}}}\neq X_{\tau_{D}}\right]=\int_{\overline{D}^{c}}\int_{D}G_{D}(x,y)J(y,z)dy f(z)dz,\quad x\in D.
\eeq
We define $K_{D}(x,z)=\int_{D}G_{D}(x,y)J(y,z)dy$ and \eqref{eqn:Levy system} can be written as
\beq\label{eqn:Levy system2}
\E_{x}\left[f(X_{\tau_{D}}), X_{\tau_{D^{-}}}\neq X_{\tau_{D}}\right]=\int_{\overline{D}^{c}}K_{D}(x,z)f(z)dz,\quad x\in D.
\eeq

Now we state the definition of harmonic functions.
\begin{defn}\label{defn:harmonic}
A function $u:D\rightarrow [0,\infty)$ is said to be
harmonic with respect to $X^{D}$ if for every open set $B$ whose closure is a compact subset of $D$,
$$
u(x)=\E_{x}\left[u(X^{D}_{\tau_{B}})\right] \quad \text{for every } x\in B.
$$
\end{defn}

Recall that an open set $D$ in $\R^{d}$ is said to be a (uniform) $C^{1,1}$ open set
if there are (localization radius) $r_{0}>0$
and $\Lambda_{0}$ such that for every $z\in\partial D$ there exist
a $C^{1,1}$ function $\phi=\phi_{z}:\R^{d}\rightarrow \R$ satisfying $\phi(0, \cdots, 0)=0$,
$\nabla\phi(0)=(0,\cdots,0)$, $|\nabla\phi(x)-\nabla\phi(y)|\leq \Lambda_{0}|x-y|$, and an orthonormal coordinate system
$CS_{z}:y=(y_{1},\cdots,y_{d-1},y_{d}):=(\tilde{y},y_{d})$ with its origin at $z$
such that $B(z,r_{0})\cap D =\{y= (\tilde{y},y_{d})\in B(0,r_{0}) \text{ in } CS_{z} : y_{d}>\phi(\tilde{y})\}$.
In this paper we will call the pair $(r_{0},\Lambda_{0})$ the characteristics of the $C^{1,1}$ open set $D$.

We state the result about the Martin boundary of a bounded $C^{1,1}$ open set $D$ with respect to $X^{D}$. Fix $x_{0}\in D$ and define
$$
M_{D}(x,y):=\frac{G_{D}(x,y)}{G_{D}(x_{0},y)},\quad x,y\in D, \,\, y\neq x, x_{0}.
$$
A positive harmonic function $f$ with respect to $X^{D}$ is called minimal if, whenever $g$ is a positive harmonic function with respect to $X^{D}$ with $g\leq f$,
one must have $f=cg$ for some positive constant $c$.
Now we recall the identification of the Martin boundary of bounded $C^{1,1}$ open sets $D$ with respect to killed processes $X^{D}$ with the Euclidean
boundary in \cite{KSV1}.
\begin{thm}\emph{(\cite[Theorem 1.5]{KSV1})}\label{thm:Martin}
Suppose that $D$ is a bounded $C^{1,1}$ open set in $\R^{d}$. For every $z\in \partial D$, there exists $M_{D}(x,z):=\displaystyle\lim_{y\rightarrow z}M_{D}(x,y)$. Furthermore, for every $z\in \partial D$, $M_{D}(\cdot,z)$ is a minimal harmonic function
with respect to $X^{D}$ and $M_{D}(\cdot,z_{1})\neq M_{D}(\cdot,z_{2})$ for $z_{1},z_{2}\in\partial D$, $z_1 \neq z_2$.
Thus the minimal Martin boundary of $D$ can be identified with the Euclidean boundary.
\end{thm}
Thus by the general theory of Martin boundary representation in \cite{KW} and Theorem \ref{thm:Martin}, we conclude that for every harmonic function
$u\geq 0$ with respect to $X^{D}$, there exists a unique finite measure $\mu$ supported on $\partial D$ such that $u(x)=\int_{\partial D}M_{D}(x,z)\mu(dz)$.
$\mu$ is called the Martin measure of $u$.
\newline

Finally we observe that
the Martin kernel $M_{D}(x,z)$ has the following two-sided estimates.
\begin{prop}
Suppose that $D$ is a bounded $C^{1,1}$ open set in $\R^{d}$, $d\geq 2$. Then there exist constants $c_{1}=c_{1}(d,D,\phi)$ and $c_{2}=c_{2}(d,D,\phi)$ such that
\beq\label{equn:Martin}
c_{1}\frac{\delta_{D}(x)}{|x-z|^{d}} \leq M_{D}(x,z) \leq c_{2}\frac{\delta_{D}(x)}{|x-z|^{d}}, \quad x\in D, \,\,z\in \partial D.
\eeq
\end{prop}
\pf
Let $$
g_{D}(x,y) :=\begin{cases}
\frac{1}{|x-y|^{d-2}}\left(1\wedge \frac{\delta_{D}(x)\delta_{D}(y)}{|x-y|^{2}}\right)&\mbox{when  } d\geq 3,\\
\log \left(1+\frac{\delta_{D}(x)\delta_{D}(y)}{|x-y|^{2}}\right) &\mbox{when } d=2.
\end{cases}
$$
Then it follows from \cite[Theorem 1.4]{KSV1} there exists $c_{1}=c_{1}(d,D,\phi)$ and $c_{2}=c_{2}(d,D,\phi)$ such that
$$
c_{1}g_{D}(x,y)\leq G_{D}(x,y) \leq c_{2}g_{D}(x,y).
$$
From Theorem \ref{thm:Martin} we immediately get the assertion of the proposition.
\qed
\section{Relative Fatou Theorem}\label{RFT}
Throughout this section we assume that $D$ is a bounded $C^{1,1}$ open set in $\R^{d}$, $d\geq 2$ with the characteristics $(r_{0},\Lambda_{0})$.
In this section we prove relative Fatou theorem for nonnegative harmonic functions $u$ and $v$ with respect to $X^{D}$.
For any finite and nonnegative measure $\mu$ supported on $\partial D$ we define
$$
M_{D}\mu(x):=\int_{\partial D}M_{D}(x,z)\mu(dz), \quad x\in D.
$$
Since $M_{D}(\cdot,z)$ is harmonic with respect to $X^{D}$ for $z\in \partial D$
(see Theorem \ref{thm:Martin}),
it is easy to see
that $M_{D}\mu(x)$ is nonnegative and harmonic with respect to $X^{D}$.

We use the following property of the surface measure $\sigma$, called \textit{Ahlfors regular condition} (see \cite[page 992]{MR}):
there exist constants $r_{1}=r_{1}(D,d)$, $C_{1}=C_{1}(D,d)$ and $C_{2}=C_{2}(D,d)$ such that for every $z\in\partial D$ and $r\leq r_{1}$
\begin{eqnarray}\label{eqn:Ahlfors}
C_{1}r^{d-1}&\leq& C_{1} \sigma(\partial D\cap (B(z,r)\setminus B(z,r/2)))\leq \sigma(\partial D \cap B(z,r))\nonumber\\
&\leq& C_{2} \sigma(\partial D\cap (B(z,r)\setminus B(z,r/2)))\leq C_{2}r^{d-1}.
\end{eqnarray}

Now we define the Stolz open set. For $Q\in\partial D$ and $\beta>1$, let
$$
A_{Q}^{\beta}=\{x\in D : \delta_{D}(x)<r_{0} \text{ and } |x-Q|<\beta\delta_{D}(x)\}.
$$
We say $x$ approaches $Q$ nontangentially if $x\rightarrow Q$ and $x\in A_{Q}^{\beta}$ for some $\beta>1$.

The next lemma is similar to \cite[Lemma 4.4]{MR}. Since we are working on $C^{1,1}$ open sets, the proof is simpler.

\begin{lemma}\label{lemma:positive}
Let $v(x)=M_{D}\nu(x)$, where $\nu$ is a finite and nonnegative measure on $\partial D$.
For $\nu$-almost every point $Q\in \partial D$, we have
$$
\liminf_{x\rightarrow Q}v(x)>0
$$
as $x\rightarrow Q$ nontangentially.
\end{lemma}
\pf
If $x\rightarrow Q$ nontangentially, there exists a constant $c>0$ such that
$$
\delta_{D}(x)\,\leq \,|x-Q| \,\leq\, c\,\delta_{D}(x).
$$
Fix $x\in D$ such that $|x-Q|<r_{1}$ and take $z\in B(Q,|x-Q|)\cap\partial D$.
Then $|x-z|\leq |x-Q|+|Q-z|\leq 2|x-Q|$ so that we obtain $M_{D}(x,z)\geq c M_{D}(x,Q)$ by \eqref{equn:Martin}.
This implies
\begin{eqnarray*}
v(x)&\geq& \int_{\partial D \cap B(Q,|x-Q|)}M_{D}(x,z)\nu(dz) \,\geq\, c\, M_{D}(x,Q) \, \nu(B(Q,|x-Q|) \cap \partial D) \\
 &\geq& c \, \frac{\delta_{D}(x)}{|x-Q|^{d}} \,\nu(B(Q,|x-Q|)\cap \partial D) \,\geq\,  c\,  \frac{\nu(B(Q,|x-Q|)\cap \partial D)}{|x-Q|^{d-1}}.
\end{eqnarray*}
By \eqref{eqn:Ahlfors} we have
$$
\frac{\sigma(B(Q,|x-Q|)\cap \partial D)}{\nu(B(Q,|x-Q|)\cap \partial D)}\,\geq\, c\, \frac{1}{v(x)},
$$
and by \cite[Theorem 5]{Be} the symmetric derivative
$$
\limsup_{x\rightarrow Q}\frac{\sigma(B(Q,|x-Q|)\cap \partial D)}{\nu(B(Q,|x-Q|)\cap \partial D)}
$$
is finite $\nu$-almost every point $Q\in\partial D$.
\qed

\begin{remark}
It follows from Lemma \ref{lemma:positive} that as $x\rightarrow Q$ nontangentially $\displaystyle\lim_{x\rightarrow Q}\frac{\delta_{D}(x)}{v(x)}=0$ for $\nu$-almost every $Q\in\partial D$, where
$v(x)=M_{D}\nu(x)$.
We will use this fact in Lemma  \ref{lemma:con1}, Theorems \ref{thm:RFT2}, and \ref{thm:RFT3}.
\end{remark}

The next lemma is an analogue of \cite[Lemma 4.3]{MR}.
\begin{lemma}\label{lemma:con1}
Let $Q\in\partial D$ and $v$ be a nonnegative harmonic function with respect to $X^{D}$ with Martin measure $\nu$.
Suppose that $\mu$ is a nonnegative finite measure on $\partial D$.
If $\displaystyle\lim_{x\rightarrow Q}\frac{\delta_{D}(x)}{v(x)}=0$, then for every $\eps>0$ we have
$$
\lim_{x\rightarrow Q}\frac{\int_{\partial D \cap \{|Q-z|\geq \eps\}}M_{D}(x,z)\mu(dz)}{v(x)}=0.
$$
If we assume
$\displaystyle\lim_{x\rightarrow Q}\frac{\delta_{D}(x)}{v(x)}=0$ nontangentially,
then the limit above must be taken nontangentially.
\end{lemma}
\pf
If $|Q-z|\geq\eps$ and $|x-Q|\leq \eps/2$, then $|x-z|\geq \eps/2$.
Thus from \eqref{equn:Martin} we have
\begin{eqnarray*}
\int_{\partial D \cap \{|Q-z|\geq\eps\}}M_{D}(x,z)\mu(dz) &\leq& c\int_{\partial D\cap\{|Q-z|\geq \eps\}}\frac{\delta_{D}(x)}{|x-z|^{d}}\mu(dz)\\
&\leq&c \, \eps^{-d}\,\delta_{D}(x)\,\mu(\partial D).
\end{eqnarray*}
Hence we have
$$
\frac{\int_{\partial D \{|Q-z|\geq \eps\}}M_{D}(x,z)\mu(dz)}{v(x)}\,\leq \, c \, \frac{\mu(\partial D)}{\eps^{d}}
\frac{\delta_{D}(x)}{v(x)}\rightarrow 0
$$
as $x\rightarrow Q$.
\qed

\begin{remark}
Note that the condition $\displaystyle\lim_{x\rightarrow Q}\frac{\delta_{D}(x)}{v(x)}=0$
cannot be omitted. To see this,
take any points $P,Q\in \partial D$ with $P\neq Q$. Let
$\mu=\nu=\delta_{\{P\}}$, $v(x)=M_{D}\nu(x)$, and $\eps=|P-Q|/2$. Then from \eqref{equn:Martin},
$\displaystyle\liminf_{x\rightarrow Q}\frac{\delta_{D}(x)}{v(x)} >0$.
Clearly
$\frac{\int_{\partial D \cap \{|Q-z|\geq \eps\}}M_{D}(x,z)\mu(dz)}{v(x)}=1$ for any $x\in D$.
\end{remark}

The next theorem is the first main result of this section.
\begin{thm}\label{thm:RFT2}
Let $Q\in\partial D$
and $v$ be a nonnegative harmonic function with respect to $X^{D}$ with Martin measure $\nu$.
Suppose that $u(x)=M_{D}\mu(x)$ for some finite nonnegative measure $\mu$ on $\partial D$ satisfying $d\mu=fd\nu$ where $f\in L^{1}(\partial D, \nu)$ and that $f$ is continuous at $Q$.
If $\displaystyle\lim_{x\rightarrow Q}\frac{\delta_{D}(x)}{v(x)}=0$, then we have
$$
\lim_{x\rightarrow Q}\frac{u(x)}{v(x)}=f(Q).
$$
If we assume
$\displaystyle\lim_{x\rightarrow Q}\frac{\delta_{D}(x)}{v(x)}=0$ nontangentially,
then the limit above must be taken nontangentially.
\end{thm}
\pf
For $Q \in \partial D$, we have
$$
\left|\frac{u(x)}{v(x)}-f(Q)\right|
\,\leq\, \frac{1}{v(x)}\int_{\partial D}|f(z)-f(Q)|M_{D}(x,z)\nu(dz)
\,\leq\,\frac{I_{1}(x)}{v(x)}+\frac{I_{2}(x)}{v(x)},
$$
where $\ds I_{1}(x):=\int_{\partial D\cap \{|z-Q|\geq \eps\}}|f(z)-f(Q)|\,M_{D}(x,z)\,\nu(dz)$ and
$\ds I_{2}(x):=\int_{\partial D\cap \{|z-Q|< \eps\}}\,|f(z)-f(Q)|\,M_{D}(x,z)\,\nu(dz)$.
Now it follows that
$\frac{I_{2}(x)}{v(x)}\leq \sup_{ \{|z-Q|<\eps\} \cap \partial D }|f(z)-f(Q)|$.
For given $\eps>0$ we can choose $r=r(\eps) >0$ such that
$|f(z)-f(Q)|<\eps$ when $|z-Q|<r$. From Lemma \ref{lemma:con1} for this $r>0$ we can take $\delta>0$ such that
if $|x-Q|<\delta$ then $\frac{I_{1}(x)}{v(x)}<\eps$ .
Thus we have $\left|\frac{u(x)}{v(x)}-f(Q)\right|<2\eps$ and this finishes the proof.
\qed

Next we consider the case when $d\mu=fd\nu+d\mu_{s}$, where $f\in L^{1}(\partial D,\nu)$ and $\mu_{s}$ is singular to $\nu$.
Consider all points $Q\in\partial D$ for which
\beq\label{eqn:Leb}
\lim_{r\rightarrow 0}\frac{\int_{B(Q,r)\cap \partial D}\left(|f(z)-f(Q)|\nu(dz)+\mu_{s}(dz)\right)}{\nu(B(Q,r)\cap\partial D)}=0.
\eeq
It is well-known that the set $Q\in \partial D$ that satisfies \eqref{eqn:Leb} is of full measure $\nu$ (for example, see \cite[Theorem 3.20 and 3.22]{Fo}).

The next lemma is the nontangential maximal inequality that is analogous to
\cite[Lemma 4.5]{MR}.
\begin{lemma}\label{lemma:maximal}
Suppose that $\mu$ and $\nu$ are nonnegative finite measures on $\partial D$.
For any $x\in D$, $Q\in \partial D$ and $t>0$ such that $|x-Q|\leq t\delta_{D}(x)$,
there exist constants $C=C(t,Q)$, $c=c(t,Q)$ such that
$$
c\,\inf_{r>0}\frac{\mu(B(Q,r)\cap\partial D)}{\nu(B(Q,r)\cap\partial D)}\,\leq\,\frac{\int_{\partial D}M_{D}(x,z)\mu(dz)}{\int_{\partial D}M_{D}(x,z)\nu(dz)}\,\leq\, C\,\sup_{r>0}\frac{\mu(B(Q,r)\cap\partial D)}{\nu(B(Q,r)\cap\partial D)}.
$$
\end{lemma}
\pf
The proof is similar to
\cite[Lemma 4.5]{MR}
but we will provide the
details for the reader's convenience.
Define $B_{n} := B(Q,2^{n}|x-Q|) \cap \partial D$ for $n\geq 1$ and $A_{n}=B_{n}\setminus B_{n-1}$ for $n\geq 2$.
Let $n_{0}$
be the smallest integer such that $2^{n_{0}}|x-Q|\geq diam(D)$. Then we have
$$
\int_{\partial D}M_{D}(x,z)\mu(dz)=\sum_{n=2}^{n_{0}}\int_{ A_{n}}M_{D}(x,z)\mu(dz) +\int_{B_{1}}M_{D}(x,z)\mu(dz).
$$
For $z\in B_{1}$ we have $|z-Q|<2|x-Q|$, which implies
$$
\frac{1}{t}|x-Q|\leq \delta_{D}(x)\leq |x-z|\leq |x-Q|+|Q-z|\leq 3|x-Q|.
$$
It follows from \eqref{equn:Martin} that
$$
\frac{M_{D}(x,z)}{M_{D}(x,Q)}\asymp \frac{|x-Q|^{d}}{|x-z|^{d}} \asymp c_{1}
$$
for some constant
$c_{1}=c_{1}(t)$.
Hence we have
$$
\inf_{z\in B_{1}}M_{D}(x,z)\leq a_{1}:=\sup_{z\in B_{1}}M_{D}(x,z)\leq c_1 \inf_{z\in B_{1}}M_{D}(x,z).
$$
It follows that
$$
\int_{B_{1}} M_{D}(x,z) \mu(dz) \,\leq\, a_{1}\,\mu(B_{1} )\,\leq\, c_1 \int_{B_{1}}M_{D}(x,z)\mu(dz).
$$
Similarly for $z\in A_{n}$ we have
$$
2^{n-2}|x-Q|\leq (2^{n-1}-1)|x-Q|\leq |z-Q|-|x-Q|\leq |x-z|\leq |x-Q|+|Q-z|\leq 2^{n+1}|x-Q|.
$$
For $z'\in A_{n}$, from \eqref{equn:Martin} we also have
$c_{2}\,M_{D}(x,z')\,\leq\, M_{D}(x,z)\,\leq\, c_{3}\,M_{D}(x,z')$. Therefore we have
$$
\inf_{z\in A_{n}}M_{D}(x,z)\leq a_{n}:=\sup_{z\in A_{n}}M_{D}(x,z)\leq c_{3}\inf_{z\in A_{n}}M_{D}(x,z).
$$
Thus we obtain
$$
a_{n}\,\mu(A_{n})\,\leq \,c_3 \,\int_{A_{n}} M_{D}(x,z)\mu(dz)
\,\leq\, c_{3} \, a_{n}\,\mu(A_{n}).
$$
Combining the estimates above we conclude that
\beq\label{eqn:max1}
c_{4}\left(\sum_{n=2}^{n_{0}}a_{n}\mu(A_{n})+a_{1}\mu(B_{1})\right) \,\leq\,
\int_{\partial D}M_{D}(x,z)\mu(dz)\,\leq\, c_{5}
\left(\sum_{n=2}^{n_{0}}a_{n}\mu(A_{n})+a_{1}\mu(B_{1})\right).
\eeq

Now define $b_{n}=\sup_{k\geq n}a_{k}$
for $n\ge 1$.
Clearly $b_{n}\geq a_{n}$ and $b_{n-1}-b_{n}\geq 0$.
Let $z'\in A_{n}$ and $z\in A_{k}$ for $k>n$. Then we have
$$
|x-z'|\leq |x-Q|+|Q-z'|\leq 2^{n+1}|x-Q|,
$$
and
$$
|x-z|\geq |z-Q|-|x-Q|\geq (2^{n}-1)|x-Q|\geq 2^{n-1}|x-Q|,
$$
so it follows that $|x-z'|\leq 4|x-z|$. From \eqref{equn:Martin} we have $M_{D}(x,z)\,\leq\, c\,M_{D}(x,z')$.
This implies that $a_{k}\,\leq \,c\,a_{n}$ for $k>n$, so $b_{n}\,\leq \,c\,a_{n}$, which in turn shows that
$\frac{b_{n}}{c}\leq a_{n}\leq b_{n}$ for $n\geq 1$.
Combining this with \eqref{eqn:max1} we have
\begin{eqnarray*}
&&c_{4}\left(\sum_{n=2}^{n_{0}}a_{n}\mu(A_{n})+a_{1}\mu(B_{1})\right)\\
&\leq&\int_{\partial D}M_{D}(x,z)\mu(dz)\\
&\leq&c_{5}\left(\sum_{n=2}^{n_{0}}a_{n}\mu(A_{n})+a_{1}\mu(B_{1})\right)\\
&\leq&c_{5}\left(\sum_{n=2}^{n_{0}}b_{n}\left(\mu(B_{n})-\mu(B_{n-1})\right)+b_{1}\mu(B_{1})\right)\\
&=&c_{5}\left(\sum_{n=2}^{n_{0}}(b_{n-1}-b_{n})\mu(B_{n-1})+b_{n_{0}}\mu(B_{n_{0}})\right),
\end{eqnarray*}
and the same estimate holds for $\nu$.
Thus we obtain
\begin{eqnarray*}
\int_{\partial D}M_{D}(x,z)\mu(dz)
&\leq& c_{5}\left(\sum_{n=2}^{n_{0}}(b_{n-1}-b_{n})\frac{\mu(B_{n-1})}{\nu(B_{n-1})}\nu(B_{n-1})
+b_{n_{0}}\frac{\mu(B_{n_{0}})}{\nu(B_{n_{0}})}\nu(B_{n_{0}})\right)\\
&\leq&c_{5}\,\sup_{r>0}\frac{\mu(B(Q,r)\cap \partial D)}{\nu(B(Q,r)\cap \partial D)}\left(\sum_{n=2}^{n_{0}}(b_{n-1}-b_{n})\nu(B_{n-1})
+b_{n_{0}}\nu(B_{n_{0}})\right)\\
&\leq&c_{6}\,\sup_{r>0}\frac{\mu(B(Q,r)\cap \partial D)}{\nu(B(Q,r)\cap \partial D)}\int_{\partial D}M(x,z)\nu(dz)
\end{eqnarray*}
and similarly we also have
$$
\int_{\partial D}M_{D}(x,z)\mu(dz) \,\geq\,  c_{7}\,\inf_{r>0}\frac{\mu(B(Q,r)\cap \partial D)}{\nu(B(Q,r)\cap \partial D)}\int_{\partial D}M(x,z)\nu(dz),
$$
and this
completes the proof.
\qed

Now we state the second main theorem of this section.
\begin{thm}\label{thm:RFT3}
Let $u,v$ be nonnegative and harmonic functions with respect to $X^{D}$.
Let $u(x)=\int_{\partial D}M_{D}(x,z)\mu(dz)$ and $v(x)=\int_{\partial D}M_{D}(x,z)\nu(dz)$, where $\mu$ and $\nu$ are nonnegative and finite measures on $\partial D$. Assume that $d\mu=fd\nu+d\mu_{s}$, where $f\in L^{1}(\partial D,\nu)$ and $\mu_{s}$ is singular to $\nu$.
Suppose $\displaystyle\lim_{x\rightarrow Q}\frac{\delta_{D}(x)}{v(x)}=0$ as $x\to Q$ nontangentially.
Then for $\nu$-almost every point $Q\in\partial D$ we have
$$
\lim_{x\rightarrow Q}\frac{u(x)}{v(x)}=f(Q)
$$
as $x\rightarrow Q$ nontangentially.
More precisely, the convergence holds for every $Q\in \partial D$ satisfying \eqref{eqn:Leb} and
$\displaystyle\lim_{x\rightarrow Q}\frac{\delta_{D}(x)}{v(x)}=0$.
\end{thm}
\pf
The proof is similar to
\cite[Theorem 4.2]{MR} but we provide the details for the reader's convenience.
Fix a point $Q\in\partial D$ that satisfies \eqref{eqn:Leb}.
Let $\eps>0$ and define $d\tilde{\mu}=|f(\cdot)-f(Q)|d\nu+d\mu_{s}$. Then we have
\begin{eqnarray*}
&&\left|\frac{u(x)}{v(x)}-f(Q)\right|\leq \frac{\int_{\partial D}M_{D}(x,z)\, \tilde{\mu}(dz)}{v(x)} \\
&=&\frac{\int_{\partial D\cap \{|z-Q|\geq \eps \}}M_{D}(x,z) \,\tilde{\mu}(dz)}{v(x)}+\frac{\int_{\partial D}M_{D}(x,z)\,\tilde{\mu}|_{B(Q,\eps)}(dz)}{v(x)},
\end{eqnarray*}
where $\tilde{\mu}|_{B(Q,\eps)}$ is the truncation of $\tilde{\mu}$ to $B(Q,\eps)\cap \partial D$.
Since $|f(\cdot)-f(Q)|\in L^{1}(\nu)$, it follows from Lemma \ref{lemma:con1} that
$$
\lim_{x\rightarrow Q}\frac{\int_{\partial D\cap \{|z-Q|\geq \eps \}}M_{D}(x,z) \,\tilde{\mu}(dz)}{v(x)}=0.
$$
Applying Lemma \ref{lemma:maximal} to the measures $\tilde{\mu}|_{B(Q,\eps)}$ and $\nu$, we get
\begin{eqnarray*}
&&\limsup_{x\rightarrow Q}\left|\frac{u(x)}{v(x)}-f(Q)\right|\\
&\leq&\limsup_{x\rightarrow Q}\frac{\int_{\partial D}M_{D}(x,z) \, \tilde{\mu}|_{B(Q,\eps)}(dz)}{v(x)}\\
&\leq&c\sup_{r>0}\frac{\tilde{\mu}|_{B(Q,\eps)}(B(Q,r)\cap \partial D)}{\nu(B(Q,r)\cap \partial D)}\\
&=&c\sup_{r\leq \eps}\frac{\int_{\partial D\cap B(Q,r)}\left(|f(z)-f(Q)|\nu(dz)+\mu_{s}(dz)\right)}{\nu(B(Q,r)\cap \partial D)}.
\end{eqnarray*}
Letting $\eps\rightarrow 0$ and using \eqref{eqn:Leb} we obtain the desired result.
\qed

In the rest of this section, we investigate relative Fatou theorem when the normalizing function corresponds
to the Martin integral with respect to the surface measure of $\partial D$.
Let $h(x):=\int_{\partial D}M_{D}(x,z)\sigma(dz)$, where $\sigma$ is the surface measure of $\partial D$.
We now prove a simple lemma about $h(x)$, which is an analogue of \cite[Equation (11)]{BD}.

\begin{lemma}\label{lemma:bdd}
There exist constants $C_{3}, C_{4}$ depending only on $D, d, \phi , x_{0}$ such that
$$
0<C_{3}\leq h(x) \leq C_{4} <\infty.
$$
\end{lemma}
\pf
For each $x\in D$ let $P=P(x)\in \partial D$ be a point such that $|x-P|=\delta_{D}(x)$.
Let $A_{n}=A_{n}(x)=\{Q\in\partial D : 2^{n-1}|x-P|\leq |x-Q| < 2^{n} |x-P|\}$.
Since $D$ is bounded, there exists $N=N(x)$
such that $\partial D\subset \bigcup_{n=1}^{N}A_{n}$.

We start with the lower bound.
By the definition of $A_n$, we get
\begin{eqnarray*}
h(x)&=&\int_{\partial D}M_{D}(x,z)\sigma(dz)\geq\int_{A_1} M_{D}(x,z)\sigma(dz)\\
&\geq&c_{3}\int_{A_1}\frac{\delta_{D}(x)}{|x-z|^{d}}\sigma(dz)\,\geq\, c_{4}\int_{A_1}\frac{\delta_{D}(x)}{2^{d}|x-P|^{d}}\sigma(dz)\\
&\geq&c_{5}2^{-d}\delta_{D}(x)|x-P|^{-d}|x-P|^{d-1}=c_{5}2^{-d}:=C_{3}.
\end{eqnarray*}

We now prove the upper bound.
From \eqref{equn:Martin} and \eqref{eqn:Ahlfors}, we have
\begin{eqnarray*}
h(x)&=&\int_{\partial D}M_{D}(x,z)\sigma(dz)\,\leq\,\sum_{n=1}^{N}\int_{A_{n}}M_{D}(x,z)\sigma(dz)\\
&\leq&c_{1}\sum_{n=1}^{N}\int_{A_{n}}\frac{\delta_{D}(x)}{|x-z|^{d}}\sigma(dz)\,\leq\, c_{1}\sum_{n=1}^{N}\int_{A_{n}}\delta_{D}(x)\left(2^{n-1}|x-P|\right)^{-d}\sigma(dz)\\
&\leq&c_{1}\,\delta_{D}(x)^{1-d}\sum_{n=1}^{N}2^{-d(n-1)}\sigma( A_n )\,\leq\, c_{1}\,\delta_{D}(x)^{1-d}\sum_{n=1}^{N}2^{-d(n-1)}\,(2^{n}|x-P|)^{d-1}\\
&\leq&c_{2}\sum_{n=1}^{N}2^{-n}\,\leq\, c_{2}\sum_{n=1}^{\infty}2^{-n}\,:=\,C_{4}.
\end{eqnarray*}
\qed

Recall that the Lebesgue set of a
function $f\in L^{1}(\partial D,\sigma)$
is the set of all $Q\in \partial D$ satisfying
$$
\lim_{r\rightarrow 0^{+}}\frac{\int_{\partial D \cap B(Q,r)}|f(z)-f(Q)|\sigma(dz)}{\sigma(\partial D\cap B(Q,r))}=0.
$$
In particular for such points $Q\in \partial D$ we have
$$
\sup_{r>0}\frac{\int_{\partial D \cap B(Q,r)}|f(z)-f(Q)|\sigma(dz)}{\sigma(\partial D\cap B(Q,r))}<\infty.
$$
It is well known that for $f\in L^{1}(\partial D,\sigma)$ the set of all Lebesgue points is of full measure $\sigma$.

For $f\in L^{1}(\partial D, \sigma)$ we define $u(x):=\int_{\partial D}M_{D}(x,z)f(z)\sigma(dz)$.
From Lemma \ref{lemma:bdd}, $\displaystyle\lim_{x\rightarrow Q}\frac{\delta_{D}(x)}{h(x)}=0$ for every $Q\in\partial D$.
Hence when the normalizing function $h(x)$ corresponds to the Martin integral with respect to the surface measure of $\partial D$, Theorem \ref{thm:RFT3} is read as the next corollary.
\begin{corollary}\label{cor:surface}
Let $f\in L^{1}(\partial D,\sigma)$ and $u(x)=\int_{\partial D}M_{D}(x,z)f(z)\sigma(dz)$.
For every Lebesgue point $Q$ of $f$,
$$
\lim_{x\rightarrow Q}\frac{u(x)}{h(x)}=f(Q), \text{ where the limit is taken nontangentially.}
$$
\end{corollary}

\section{Fatou theorem in a ball}\label{ball}

In this section we prove that Fatou theorem holds for $X$ when $D$ corresponds to a ball centered at $x_{0}$. That is, we prove that
for any nonnegative harmonic functions $u(x)$ with respect to $X^{B(x_{0},r)}$ on a ball $B(x_{0},r)$, the nontangential limit
$$
\lim_{A_{Q}^{\beta}\ni x\rightarrow Q\in\partial B(x_{0},r)}u(x) \text{ exists for } \sigma\text{-a.e. } Q\in \partial B(x_{0},r).
$$

In order to prove Fatou theorem for a ball, we first establish a few lemmas.
Take a (smooth) open sets $D_{n}$ such that $D_{n}\subset \overline{D_{n}}\subset D_{n+1}$ and $\cup_{n=1}^{\infty}D_{n}=D$. Let $\tau_{n}=\tau_{D_{n}}$.

\begin{lemma}\label{lemma:exit on the boundary}
For each $x\in D$ we have
$$
\P_{x}\left(\{w\in \Omega : \cap_{n=1}^{\infty}\{\tau_{n}<\tau_{D}\} =\{X_{\tau_{D}}\in\partial D\} \}\right)=1.
$$
\end{lemma}
\pf
Since $D_{n}\subset D$ we have $\tau_{n}\leq \tau_{D}$. Also it follows from the definition of $\tau_{D}$ we have $X_{\tau_{D}}\notin D$.
Suppose that $w\in \{w\in \Omega : X_{\tau_{D}}\in \partial D\} \cap \{w\in\Omega : \tau_{n}=\tau_{D}\}$. Then
$X_{\tau_{n}}(\omega)\in\partial D$. Since $D=B(x_{0},r)$ has zero $d$-dimensional Lebesgue measure it follows from \eqref{eqn:Levy system2} we have
$$
\P_{x}\left(X_{\tau_{n}}\in \partial D\right)=\P_{x}\left(X_{\tau_{n}}\in \partial D, X_{\tau_{n}^{-}}\neq X_{\tau_{n}}\right)=\int_{\partial D} K_{D_{n}}(x,z)dz=0.
$$
Hence $\{w\in \Omega : X_{\tau_{D}}\in \partial D\}\subset \{w\in\Omega : \tau_{n}<\tau_{D}\}$ $\P_{x}$-a.e. $x\in D$ for all $n\in \mathbb{N}$.
Now take $w\in \cap_{n=1}^{\infty}\{w\in\Omega : \tau_{n}<\tau_{D}\}$ then we have $X_{\tau_{D}^{-}}(\omega)=X_{\tau_{D}}(\omega)\in \partial D$.
\qed

For any $Q\in\partial D$, we let $\phi_{Q}$ be the $C^{1,1}$ function associated with $Q$ in the definition
of $C^{1,1}$ open sets. For any $x\in \{y=(y,y_{d})\in B(Q,R): y_{d}>\phi_{Q}(\tilde{y})\}$ we put
$\rho_{Q}(x):=x_{d}-\phi_{Q}(\tilde{x})$.
For $r,s>0$, we define
$$
D_{Q}(r,s):=\{y\in D : r>\rho_{Q}(y)>0, |\tilde{y}|<s\}.
$$
Let $r_{1}:=r_{0}/4(\sqrt{1+(1+\Lambda_{0})^{2}})$.
The following result is \cite[Lemma 4.3]{KSV1}.
\begin{lemma}\emph{(\cite[Lemma 4.3]{KSV1})}\label{lemma:exit inside}
There exist constants $\lambda_{0}>2r_{1}^{-1}$, $\kappa_{0}\in(0,1)$,
and $c=c(r_{0},\Lambda_{0})$ such that for every
$\lam\geq\lambda_{0}$, $Q\in\partial D$, and $x\in D_{Q}(2^{-1}(1+\Lambda_{0})^{-1}
\kappa_{0}\lam^{-1},\kappa_{0}\lam^{-1})$
with $\tilde{x}=0$,
$$
\P_{x}\left(X_{\tau_{D_{Q}(\kappa_{0}\lambda^{-1},\lambda^{-1})}}\in D\right)\leq c\lambda\delta_{D}(x).
$$
\end{lemma}

Next we prove that as the starting point $x$ approaches $Q\in \partial D$ the probability of exiting $D$ near the point $Q$ under the condition the processes exit $D$ through the boundary will converge to 1.
\begin{lemma}\label{lemma:exit}
For any $r<r_{0}$, $Q\in\partial D$, and $\eps>0$ there exists a constant $r_{2}>0$
such that for any $x\in D$ with $|x-Q|<r_{2}$
we have
$$
\P_{x}\left(X_{\tau_{D}}\notin B(Q,r), X_{\tau_{D}}\in \partial D\right)<\eps.
$$
\end{lemma}
\pf
For any given $r<r_{0}$, we take a large enough $\lambda$ so that
$D_{Q}(\kappa_{0}\lambda^{-1},\lambda^{-1})\subset B(Q,r)$.
Then
$$
\{X_{\tau_{D}}\notin B(Q,r), X_{\tau_{D}}\in \partial D\}\subset \{X_{\tau_{D_{Q}(\kappa_{0}\lambda^{-1},\lambda^{-1})}} \in D\}.
$$
It follows from Lemma \ref{lemma:exit inside} that
$$
\P_{x}\left(X_{\tau_{D}}\notin B(Q,r), X_{\tau_{D}}\in \partial D\right)\leq \P_{x}\left(X_{\tau_{D_{Q}(\kappa_{0}\lambda^{-1},\lambda^{-1})}} \in D\right)\leq c\lambda \delta_{D}(x)\leq c\lambda |x-Q|.
$$
Taking $r_{2}=\eps/c\lambda$, we arrive at the desired assertion.
\qed

Recall that the point $x_{0}\in D$ is the point such that $M_{D}(x_{0},y)=1$ for all $y\in D$ (hence $M_{D}(x_{0},z)=1$ for all $z\in\partial D$ as well).
We prove that the Martin kernel $M_{D}(x,z)$ is the Radon-Nikodym derivative with respect to harmonic measures supported on the boundary of $D$ with different starting points.

\begin{lemma}\label{lemma:Martin kernel as RN derivative}
Let $x\in D$. Then for each Borel set $A\subset \partial D$ we have
$$
\P_{x}\left(X_{\tau_{D}}\in A\right)=\int_{A}M_{D}(x,z)\P_{x_{0}}(X_{\tau_{D}}\in dz).
$$
\end{lemma}
\pf
For each Borel set $A\subset\partial D$ the function $x\rightarrow \P_{x}\left(X_{\tau_{D}}\in A\right)$ is harmonic with respect to $X^{D}$. Hence it follows from Theorem \ref{thm:Martin} there exists a nonnegative finite measure $\mu_{A}$ supported on $\partial D$ such that
$$
\P_{x}\left(X_{\tau_{D}}\in A\right)=\int_{\partial D}M_{D}(x,z)\mu_{A}(dz), \quad x\in D.
$$
For Borel sets $A, B\subset \partial D$, $A\cap B=\phi$ we have $\P\left(X_{\tau_{D}}\in A\cup B\right)=\int_{\partial D}M_{D}(x,z)\mu_{A\cup B}(dz)$ and $\P\left(X_{\tau_{D}}\in A\cup B\right)=\P\left(X_{\tau_{D}}\in B\right)+\P\left(X_{\tau_{D}}\in  B\right)
=\int_{\partial D}M_{D}(x,z)(\mu_{A}+\mu_{B})(dz)$. By the uniqueness of the Martin representation we have
$$
\mu_{A\cup B}=\mu_{A}+\mu_{B}, \quad \text{if } A\cap B=\phi.
$$

We will prove that $\text{supp}(\mu_{A})\subset \bar{A}$.
For each $\eps>0$ define $\mu_{A}^{\eps}=\mu_{A}|_{A^{\eps}}$, where $A^{\eps}=\{z\in\partial D : \text{dist}(z,A)>\eps\}$.
Let $f(x):=\int_{\partial D}M_{D}(x,z)\mu_{A}^{\eps}(dz)$.
For $Q\in \partial D \setminus A^{\eps/2}$ it follows from \eqref{equn:Martin} and the dominated convergence theorem
$$
\lim_{x\rightarrow Q\in \partial D\setminus A^{\eps/2}}f(x)=0.
$$
For $w\in A^{\eps/2}$ it follows from Lemma \ref{lemma:exit} and the fact $f(x)\leq \P_{x}(X_{\tau_{D}}\in A)$ we have
$$
\lim_{x\rightarrow w\in A^{\eps/2}}f(x)\leq \lim_{x\rightarrow w\in A^{\eps/2}}\P_{x}(X_{\tau_{D}}\in A)\leq \lim_{x\rightarrow w\in A^{\eps/2}}\P_{x}\left(X_{\tau_{D}}\notin B(w,\eps/2), X_{\tau_{D}}\in \partial D\right)=0.
$$
Hence we have proved that
\beq\label{eqn:f=0 on the boundary}
\lim_{x\rightarrow Q\in \partial D}f(x)=0 \text{ for all } Q\in\partial D.
\eeq
Now we will show that $f(x)=0$ for all $x\in D$ and this will imply that $\mu_{A}^{\eps}=0$.
Take $D_{n}$ as before. Since $f(x)$ is harmonic with respect to $X^{D}$ we have
\begin{eqnarray}\label{eqn:MK limit}
f(x)&=&\E_{x}[f(X^{D}_{\tau_{D_{n}}})]=\E_{x}[f(X_{\tau_{n}}), \tau_{n}<\tau_{D}]\nonumber\\
&=&\E_{x}[\int_{\partial D}M_{D}(X_{\tau_{n}},z)\mu_{A}^{\eps}(dz), \tau_{n}<\tau_{D}].
\end{eqnarray}
Since $\E_{x}[\int_{\partial D}M_{D}(X_{\tau_{n}},z)\mu_{A}^{\eps}(dz), \tau_{n}<\tau_{D}]=f(x)\leq \P_{x}(X_{\tau_{D}}\in A)$ for all $n\in \mathbb{N}$, it follows from Lebesgue dominated convergence theorem, Lemma \ref{lemma:exit on the boundary}, and \eqref{eqn:f=0 on the boundary} we have
\begin{eqnarray*}
&&\lim_{n\rightarrow \infty}\E_{x}[\int_{\partial D}M_{D}(X_{\tau_{n}},z)\mu_{A}^{\eps}(dz), \tau_{n}<\tau_{D}]\\
&=&\E_{x}[\lim_{n\rightarrow \infty}\int_{\partial D}M_{D}(X_{\tau_{n}},z)\mu_{A}^{\eps}(dz), X_{\tau_{D}}\in \partial D ]\\
&=&0.
\end{eqnarray*}
This shows that $f(x)=0$ and $\mu_{A}^{\eps}=0$. Since $\eps>0$ is arbitrary we have $\text{supp}(\mu_{A})\subset\bar{A}$.
Now the rest of the proof is identical to \cite[Lemma 3.1]{Kim2} and will be omitted.
\qed

In the rest of the section, we let $D=B(x_{0},r)$ for some $0<r<\infty$.
Since $X$ is rotationally invariant, $\P_{x_{0}}(X_{\tau_{B(x_{0},r)}}\in dy, X_{\tau_{B(x_{0},r)}}\in \partial B(x_{0},r))$ must be a normalized surface measure $S^{-1}\sigma(dz)$, where $S=S(X,D):=\P_{x_{0}}(X_{\tau_{B(x_{0},r)}}\in \partial B(x_{0},r))$.
Hence for any Borel set $A\subset\partial D$, it follows from Lemma \ref{lemma:Martin kernel as RN derivative} we have
$$
\P_{x}\left(X_{\tau_{B(x_{0},r)}}\in A\right)=\int_{\partial B(x_{0},r)}M_{B(x_{0},r)}(x,z)S^{-1}\sigma(dz).
$$
It follows from \cite[Theorem 3.2]{Mi} and the remark below that for each $Q\in \partial B(x_{0},r)$ we have
\beq\label{eqn:limit=1}
\lim_{x\rightarrow Q\in \partial B(x_{0},r)}\P_{x}(X_{\tau_{B(x_{0},r)}}\in \partial B(x_{0},r))=1.
\eeq
Let $u(x)=\int_{\partial B(x_{0},r)}M_{B(x_{0},r)}(x,z)f(z)\sigma(dz)$, $f(z)\in L^{1}(\partial B(x_{0},r), \sigma)$ and $v(x)=\int_{\partial B(x_{0},r)}M_{B(x_{0},r)}(x,z)\sigma(dz)$.
Then clearly $u(x)$ and $v(x)$ are harmonic with respect to $X^{B(x_{0},r)}$. Also it follows from Lemma \ref{lemma:Martin kernel as RN derivative}
$v(x)=S\P_{x}(X_{\tau_{B(x_{0},r)}}\in\partial B(x_{0},r))$ and $\lim_{x\rightarrow Q\in\partial D}v(x)=S$ for all $Q\in\partial D$.
Now we prove Fatou theorem for a ball.
\begin{thm}
Let $u(x)=\int_{\partial B(x_{0},r)}M_{B(x_{0},r)}(x,z)f(z)\sigma(dz)$, $f(z)\in L^{1}(\partial B(x_{0},r), \sigma)$. Then $\sigma$-a.e. $Q\in \partial B(x_{0},r)$ the nontangential limit exists and is equals to $Sf(Q)$. That is,
$$
\lim_{A_{Q}^{\beta}\ni x\rightarrow Q\in \partial B(x_{0},r)}u(x)=Sf(Q), \quad S=\P_{x_{0}}(X_{\tau_{B(x_{0},r)}}\in \partial B(x_{0},r)).
$$
\end{thm}
\pf
It follows from Corollary \ref{cor:surface} $\lim_{A_{Q}^{\beta}\ni x\rightarrow Q\in \partial B(x_{0},r)}\frac{u(x)}{v(x)}=f(Q)$ and from \eqref{eqn:limit=1}
$\lim_{x\rightarrow Q\in \partial B(x_{0},r)}v(x)=S$. Hence we have
$$
\lim_{A_{Q}^{\beta}\ni x\rightarrow Q\in \partial B(x_{0},r)}u(x)=\lim_{A_{Q}^{\beta}\ni x\rightarrow Q\in \partial B(x_{0},r)}\frac{u(x)}{v(x)}\cdot\lim_{A_{Q}^{\beta}\ni x\rightarrow Q\in \partial B(x_{0},r)}v(x)=Sf(Q).
$$
\qed

In \cite{Li} it is proved that there exists a bounded (classical) harmonic function on the unit disk in $\R^{2}$ that fails to have tangential limits for a.e. $\theta\in[0,2\pi]$.
Using the similar method, in \cite{Kim2, Kim} the author showed that the Stolz open sets are best possible sets for Fatou theorem and relative Fatou theorem for
transient censored stable processes and stable processes, respectively for $d=2$ and $D=B(0,1)$.

A curve $C_{0}$ is called a tangential curve in $B(0,1)$ if $C_{0}\cap \partial B(0,1)=\{w\}\in \partial B(0,1)$, $C_{0}\setminus \{w\}\subset B(0,1)$, and for any $r>0$ and $\beta>1$ $C_{0}\cap B(w,r)\nsubseteq A_{w}^{\beta}\cap B(w,r)$.
Let $C_{\theta}$ be a rotated curve $C_{0}$ about the origin through an angle $\theta$.
We will adapt arguments in \cite{Kim2, Kim, Li} to prove that the Stolz open sets are best possible sets for Fatou theorem for $X$ by showing that there exists bounded harmonic function $u(x)$ with respect to $X^{B(0,1)}$ such that the tangential limit $\lim_{x\in C_{\theta}, x\rightarrow Q}u(x)$ does not exist, where $C_{\theta}$ is a tangential curve inside $B(0,1)$.

We start with a simple lemma that is analogue to \cite[Lemma 2]{Li} (see also \cite[Lemma 3.19]{Kim2} and \cite[Lemma 3.22]{Kim}).
Let $D=B(0,1)\in \R^{2}$, $x_{0}=0$, and $\sigma_{1}$ be the normalized surface measure of $\partial B(0,1)$.
\begin{lemma}\label{lemma:near 1}
Let $h(x):=\int_{\partial B(0,1)}M_{B(0,1)}(x,z)\sigma_{1}(dz)$ and $U(z)$ be nonnegative, measurable on $\partial B(0,1)$, and $0\leq U(e^{i\theta})\leq 1$, $\theta\in [0,2\pi]$.
Suppose that $U(e^{i\theta})=1$ for $\theta_{0}-\lambda\leq\theta\leq \theta_{0}+\lambda$ for some $0<\lambda<\pi$.
Let $u(x)=\int_{\partial B(0,1)}M_{B(0,1)}(x,z)U(z)\sigma_{1}(dz)$, $x\in B(0,1)$.
Then for any $\eps>0$ there exists
$\delta=\delta(\eps,\phi)$, independent of $\lambda$, such that
$$
1-\eps \leq \frac{u(\rho e^{i\theta_{0}})}{h(\rho e^{i\theta_{0}})}\leq 1, \quad \text{if } \rho>1-\lambda\delta.
$$
\end{lemma}
\pf
Since $0\leq U(z)\leq 1$ we have
$$
0\leq \frac{u(x)}{h(x)}=\frac{1}{h(x)}\int_{\partial B(0,1)}M_{B(0,1)}(x,z)U(z)\sigma_{1}(dz)
\leq \frac{1}{h(z)}\int_{\partial B(0,1)}M_{B(0,1)}(x,z)\sigma_{1}(dz)=1.
$$
Let $V(z):=\frac{1-U(z)}{2}$ so that $0\leq V(z)\leq \frac12$ and $V(e^{i\theta})=0$ for $\theta_{0}-\lambda\leq\theta\leq \theta_{0}+\lambda$.
By the triangular inequality we have $|e^{i\theta_{0}}-e^{i\theta}|\leq |e^{i\theta_{0}}-\rho e^{i\theta_{0}}|+
|\rho e^{i\theta_{0}}- e^{i\theta}|= (1-\rho)+|\rho e^{i\theta_{0}}-e^{i\theta}|$. Hence
\begin{eqnarray*}
|\rho e^{i\theta_{0}}- e^{i\theta}|&\geq& |e^{i\theta_{0}}-e^{i\theta}|-(1-\rho)\geq 2\left|\sin(\frac{\theta_{0}-\theta}{2})\right|-\delta|\theta_{0}-\theta|\\
&\geq&\frac{2}{\pi}|\theta_{0}-\theta|-\delta|\theta_{0}-\theta|\\
&=&(\frac{2}{\pi}-\delta)|\theta_{0}-\theta|
\end{eqnarray*}
for $|\theta_{0}-\theta|>\lambda$.
Hence from \eqref{equn:Martin} we have for $\rho>1-\lambda\delta$
\begin{eqnarray*}
&&\int_{0}^{2\pi}M_{B(0,1)}(\rho e^{i\theta_{0}},e^{i\theta})V(e^{i\theta})d\theta\\
&\leq&c_{1}(1-\rho)\int_{0}^{2\pi}\frac{V(e^{i\theta})}{|\rho e^{i\theta_{0}}-e^{i\theta}|^{2}}d\theta\\
&\leq&c_{1}(1-\rho)(\frac{2}{\pi}-\delta)^{-2}\int_{|\theta-\theta_{0}|>\lambda}\frac{d\theta}{|\theta_{0}-\theta|^{2}}\\
&\leq&c_{1}\frac{1-\rho}{\lambda}(\frac{2}{\pi}-\delta)^{-2}\\
&\leq&c_{1}\frac{\delta}{(\frac{2}{\pi}-\delta)^{2}}.
\end{eqnarray*}
From Lemma \ref{lemma:bdd} $h(x)>c_{2}$ for some constant $c_{2}>0$.
Hence if $\delta\leq\frac{1}{\pi}$ we have
\begin{eqnarray*}
\frac{u(\rho e^{i\theta_{0}})}{h(\rho e^{i\theta_{0}})}&=&\frac{1}{h(\rho e^{i\theta_{0}})}\frac{1}{2\pi}\int_{0}^{2\pi}M_{B(0,1)}(\rho e^{i\theta_{0}},e^{i\theta})(1+2V(e^{i\theta}))d\theta\\
&\geq&\frac{1}{h(\rho e^{i\theta_{0}})}(h(\rho e^{i\theta_{0}})-2c_{1}\frac{\delta}{(\frac{2}{\pi}-\delta)^{2}})\\
&\geq&1-c_{3}\delta.
\end{eqnarray*}
Now for given $\eps$ take $\delta=\frac{1}{\pi}\wedge \frac{\eps}{c_{3}}$ and we reach the conclusion of the lemma.
\qed

Once we have Lemma \ref{lemma:near 1} by adapting the argument in \cite{Li} we have the following theorem.

\begin{thm}
There exists a bounded and nonnegative harmonic function $u(x)$ with respect to
$X^{B(0,1)}$ such that for a.e. $\theta\in [0,2\pi]$
$$
\lim_{|x|\rightarrow 1, x\in C_{\theta}}u(x) \text{ does not exist}.
$$
\end{thm}
\pf
Let $h(x):=\int_{\partial B(0,1)}M_{B(0,1)}(x,z)\sigma_{1}(dz)$ and $u(x)=\int_{\partial B(0,1)}M_{B(0,1)}(x,z)U(z)\sigma_{1}(dz)$, $U(z)\geq 0$.
From Corollary \ref{cor:surface} the nontangential limit $\lim_{A_{Q}^{\beta}\ni x\rightarrow Q\in\partial B(0,1)}\frac{u(x)}{h(x)}$ exists for
a.e. $Q\in\partial D$. By following the argument in \cite{Li} there exist nonnegative harmonic functions $u_{k}(x)$ with respect to $X^{B(0,1)}$
defined on some $E_{k}^{*}$ such that
$$
\lim_{x\rightarrow z\in \partial B(0,1)}\frac{u_{k}(x)}{h(x)}=0 \text{ radially and } \limsup_{x\rightarrow z\in \partial B(0,1)}\frac{u_{k}(x)}{h(x)}=2^{-k} \text{ along one branch of } C_{\theta}.
$$
Let $u(x)=\sum_{k=1}^{\infty}u_{k}(x)$. For this $u(x)$ by following the argument in \cite{Li} with Lemma \ref{lemma:near 1}
(see also \cite[Theorem 3.23]{Kim}) we have
$$
\lim_{|x|\rightarrow 1, x\in C_{\theta}}\frac{u(x)}{h(x)} \text{ does not exist for a.e. } \theta\in [0,2\pi].
$$
From Lemma \ref{lemma:Martin kernel as RN derivative} $h(x)=\frac{S}{2\pi}\P_{x}\left(X_{\tau_{B(0,1)}}\in \partial B(0,1)\right)$ and it follows from
\cite[Theorem 3.2]{Mi} the (ordinary) limit of $h(x)$ exists and $\lim_{B(0,1)\ni x\rightarrow Q\in \partial B(0,1)}h(x)=\frac{S}{2\pi}$ for all $Q\in\partial B(0,1)$.
Hence we have
$$
\lim_{|x|\rightarrow 1, x\in C_{\theta}}u(x)  \text{ does not exist for a.e. } \theta\in [0,2\pi].
$$
\qed

\medskip
{\bf Acknowledgement}
Y. Lee wants to express her heartfelt thanks to Universit\"at Bielefeld for their hospitality, where this paper was written in part.
The authors are grateful to Professor P. Kim and R. Song for many helpful comments and suggestions on this paper.

\bibliographystyle{plain}
\bibliography{RFT_ref}

\end{doublespace}

\vskip 0.3truein

{\bf Yunju Lee}

PDE and Functional Analysis Research Center, Seoul National University, Republic of Korea

E-mail: \texttt{yunjulee333@gmail.com}

\bigskip

{\bf Hyunchul Park}

Department of Mathematics, The College of William and Mary, VA 23187, USA

E-mail: \texttt{hpark02@wm.edu}

\bigskip

\end{document}